\documentclass{amsart}%{amsproc}

\newtheorem{theorem}{Theorem}[section]
\newtheorem{lemma}[theorem]{Lemma}
\newtheorem{proposition}[theorem]{Proposition}

\theoremstyle{definition}
\newtheorem{definition}[theorem]{Definition}

\newtheorem{corollary}[theorem]{Corollary}
\newtheorem{problem}[theorem]{Problem}

\theoremstyle{remark}

\newcommand{\selabel}[1]{\label{se:#1}}
\newcommand{\seref}[1]{Section~\ref{se:#1}}
\newcommand{\lelabel}[1]{\label{le:#1}}

\newcommand{\prlabel}[1]{\label{pr:#1}}
\newcommand{\prref}[1]{Proposition~\ref{pr:#1}}
\newcommand{\colabel}[1]{\label{co:#1}}
\newcommand{\coref}[1]{Corollary~\ref{co:#1}}

\newcommand{\delabel}[1]{\label{de:#1}}
\newcommand{\deref}[1]{Definition~\ref{de:#1}}
\newcommand{\eqlabel}[1]{\label{eq:#1}}
\newcommand{\equref}[1]{(\ref{eq:#1})}

\newcommand{\id}{\mbox{id}}

\newcommand{\End}{\mbox{End}}
\def\Hom{\rm Hom}

\def\ol{\overline}

\def\dul#1{\underline{\underline{#1}}}
\def\Nat{\dul{\rm Nat}}
\def\ot{\otimes}

\input diagrams

\title[Are Biseparable Extensions Frobenius?]
{Are Biseparable Extensions Frobenius?}

\author{Stefaan Caenepeel}
\address{Faculty of Applied Sciences\\
Free University of Brussels (VUB)\\
Pleinlaan 2\\
B-1050 Brussels\\ Belgium}
\email{scaenepe@vub.ac.be}
\author{Lars Kadison}
\address{Matematiska Institutionen \\ G{\" o}teborg
University \\
S-412 96 G{\" o}teborg\\ Sweden}
\email{kadison@math.ntnu.no}
\thanks{The authors thank G. Bergman,
O. Kerner, O.A. Laudal, V. Mazurchuk, D. Nikshych,
P. Seibt,
  A.A. Stolin, and R. Wisbauer for
discussions related to this paper, as well as the referees
for valuable comments.  The second author thanks the first and V.U.B. for
a nice visit in Brussels in September 2000, and NorFA for financial support.}

\subjclass{16L60, 16H05}
\date{}

\begin{document}

\begin{abstract}
In Secion~1 we describe what is known of the extent to which a separable
extension of
unital associative rings
is a Frobenius extension.  A problem of this kind is suggested by asking
if three algebraic axioms for finite Jones index subfactors are dependent.
In Section~2  the problem in the title is formulated in terms of separable
bimodules.
In Section~3 we specialize the problem to
ring extensions, noting that a biseparable extension
is a two-sided  finitely generated projective, split, separable extension.
Some reductions of the problem are discussed
and solutions in special cases are provided.
In Section~4 various examples are provided of projective separable
extensions that are neither finitely generated
nor Frobenius and which give obstructions to
weakening the hypotheses of the question in the title.
We show in Section~5 that
characterizations of the separable extensions among the Frobenius extensions in
\cite{HS,K,K99} are special cases of a result for adjoint functors.
\end{abstract}

\maketitle

%%%%%%%%%%%%%%%%%%%%%%%%%%%%%%%%%%%%%%%%%%%%%%%%%%%%%%%%%%%%%%%%%%%%%%%%%
%%%%%%%  Introduction                      %%%%%%%%%%%%%%%%%%%%%%%%%%%%%%
%%%%%%%%%%%%%%%%%%%%%%%%%%%%%%%%%%%%%%%%%%%%%%%%%%%%%%%%%%%%%%%%%%%%%%%%%
\begin{section}{Introduction}\selabel{1}

An old problem is the extent to which separable algebras
are Frobenius algebras.  By a Frobenius algebra we mean a finite dimensional
algebra $A$ with a non-degenerate linear
functional, which  induces an $A$-module isomorphism $A \cong A^*$;
symmetric algebra if this isomorphism is an $A$-bimodule map.
Eilenberg and Nakayama observed in \cite{EN} that
the (reduced) trace of a central simple
algebra over a field is non-degenerate,
which implies that a finite dimensional
semisimple algebra is  symmetric. Passing
to a commutative ground ring $k$,
Hattori \cite{H} and DeMeyer \cite{D} showed that a
$k$-projective separable $k$-algebra $A$
is symmetric as well if the Hattori-Stallings rank of $A$
over its center $C$ is an invertible element in $C$.
Endo and Watanabe  extended this result  to  $k$-projective
separable faithful $k$-algebras essentially by  using the Auslander-Goldman
Galois theory for commutative rings to define a more general
notion of reduced trace \cite{EW}.

  The main theorem in \cite{EW} led to several general results by Sugano
\cite{S70}
for when separable extensions \cite{HS} are Frobenius extensions \cite{K60}.
These are noncommutative ring extensions and are natural objects for study
from the point of view of induced representations \cite{Hoch}.
Sugano shows that a centrally projective separable extension $R/S$ is
Frobenius since it satisfies
$R \cong S \otimes_{Z(S)} C_R(S)$ where the centralizer $C_R(S)$
is faithfully
projective and separable over the center $Z(S)$,
whence Frobenius. Somewhat similarly, it
is shown that a split one-sided
finite projective H-separable extension $R/S$ is Frobenius,
since in this case the endomorphism
ring $\End (R_S) \cong R \otimes_{Z(R)} C_R(S)$ with $C_R(S)$
again separable,
the result following from
the endomorphism ring theorem as developed in \cite{K60,Mu64,M67}:
if $R_S$ is generator module, $R/S$ is   Frobenius
iff $\End(R_S) / R$ is Frobenius.  However, it is
implicit in the literature that there are several cautionary examples
showing separable extensions are not always Frobenius extensions
in the ordinary  untwisted sense \cite{K99}:
in Section~4  we show that a non-finite
ring leads to an example of split, separable, two-sided projective
extension which is not finitely generated, whence not Frobenius.

As an independent line of inquiry, algebraic axioms for finite Jones index
subfactors
have been investigated in
\cite{K95,K96,NK}.
If we simplify the discussion somewhat, we may start with
an irreducible subfactor $M/N$ of finite index:  from the Pimsner-Popa
orthonormal
base, the natural modules $M_N$ and ${}_NM$ are finite projective  \cite{JS},
and the algebra extension $M/N$ is
(1)  split,
(2) separable, and
(3)  Frobenius.
A ring extension $R/S$ is said to be \textit{split}
if there is a bimodule projection $E: R \to S$.  At the same time,
Axiom~(2) yields a Casimir element $e = \sum_i x_i \otimes y_i$
such that $\sum_i x_i y_i = 1$. A problem in the independence of the
axioms above becomes whether Axioms~(1) and (2) imply Axiom~(3)
in the presence of the assumption of two-sided finite projectivity
of the ring extension: i.e., whether
a bimodule map $E: A \to S$ and Casimir element $e \in A \otimes_S A$ may
be chosen such that
$\sum_i E(x_i)y_i = 1 = \sum_i x_iE(y_i)$. Equivalently,
can a bimodule map $E$ be found such that
the $E$-multiplication on $R \otimes_S R$ \cite{J}
is unital? That Axiom~(3) in combination with (1) or (2) does not imply the
other are easy examples disposed of in Section~3.
Many other
algebraic examples of split separable Frobenius can be found \cite{K95}
but none thus far that are 
finite projective split separable and not Frobenius.

In this paper, we will formulate the problem of independence of axioms
for subfactors in several  algebraic ways.  In Section~2
we first formulate the problem  using  separable bimodules \cite{S71},
a theory in which
separable extension and split extension become dual notions \cite{K96,K99}.
In analogy with ``bialgebra,'' we will baptise finite projective split
separable extensions
as \textit{biseparable extensions}.  Posed in the negative, our question
then becomes
if biseparable extensions are Frobenius.  This question will be formulated
in several other ways in Section~3, with one special case being answered
in the affirmative.  We point out here
that the problem has many interesting sub-problems if restrictions
are placed on the rings (e.g.,
``finite dimensional algebras,''  ``Hopf algebras,'' etc.).

In Section~5 we discuss a type of converse
to the considerations above.  We find a  common feature of the theorems in
\cite{HS,K,K99} on when a Frobenius extension or bimodule
is separable: in each case, it is a specific example of a known
theorem on adjoint functors, which we expose in this last section.
\subsection*{Preliminaries}  A ring $R$ will mean a unital associative ring.
A ring homomorphism sends 1 into 1.  A right module
$M_R$ or left module ${}_RM$ is always unitary.  Bimodules are associative with
respect to the left and right actions.

An $R$-$S$-bimodule $M$ is denoted by ${}_RM_S$.
Its right dual is defined by $M^* := \Hom(M_S,S_S)$, an
$S$-$R$-bimodule where  $sfr(m) := sf(rm)$.
The left dual of $M$ is ${}^*M := \Hom({}_RM,{}_RR)$
is also $S$-$R$-bimodule where  $(m)(sfr) := [(ms)f]r$.
Both $M \mapsto M^*$ and $M \mapsto {}^*M$ are contravariant functors of
bimodule categories, sending ${}_R{\mathcal{M}}_S \rightarrow
{}_S{\mathcal{M}}_R$.

If $R = S$ in the last paragraph, denote $\hat{M} := \Hom_{S-S}(M,S)$.
Define the group of $S$-central or Casimir elements by  \(
M^S := \{ m \in M|\ ms = sm,\ \forall s \in S \}
\).
Note that $(M^*)^S = ({}^*M)^S = \hat{M}$.

If ${}_RM_S$, ${}_RN_T$, ${}_TQ_S$ and ${}_SP_T$ are bimodules, then
$M \otimes_S P$ receives the natural $R$-$T$-bimodule structure indicated
by $r(m \otimes n)t := rm \otimes nt$,
and the group of
right module homomorphisms $\Hom(M_S,Q_S)$
receives the natural $T$-$R$-bimodule structure indicated by
$(tfr)(m) := t(f(rm))$. The group of left module homomorphisms
$\Hom({}_RM,{}_RN)$ receives the natural $S$-$T$-bimodule structure
indicated by $(m)(sft) = ((ms)f)t$.
All bimodules arising from Hom
and tensor in this paper are the natural ones unless otherwise indicated.

A {\it ring extension} $R/S$ is a ring homomorphism
$ S \stackrel{\iota}{\rightarrow} R$.
A ring extension is an {\it algebra} if $S$ is commutative and $\iota$
factors into $S \rightarrow Z(R) \hookrightarrow R$ where
$Z(R) := R^R$ is the center of $R$.
A ring extension is {\it proper}
if $\iota$ is 1-1, in which case identification is made.

The natural bimodule ${}_SR_S$ is given by $ s \cdot a \cdot s' := \iota(s)a
\iota(s')$.  In particular, we  consider the natural modules
$R_S$ and ${}_SR$.
An adjective, such as right projective or projective, for  the
ring extension $R/S$ refers to the same adjective for
one or both of these natural modules.
The structure map $\iota$ is usually  suppressed.

Separable extensions are studied in \cite{HS,K95,K,RAFAEL,S70}
among others. A ring extension $R/S$ is \textit{separable}
if the natural (multiplication) map $R \otimes_S R \to R$ is a split
epimorphism of $R$-bimodules.  Examples are abundant among finite dimensional
algebras since a separable algebra is a separable extension of any of
its subalgebras.
The next proposition, whose proof follows
Sugano \cite[Prop.\ 1]{S82},  is  important to keep in mind
when finding examples of separable extensions
from the class of finite dimensional algebras.
% because of the Wedderburn Principal Theorem.

\begin{proposition}
If $ 0 \to J \to A \stackrel{\pi}{\to} S \to 0$
is a split exact sequence of algebras
and  $A/S$ is  a separable extension,
then  $J$ is an idempotent ideal (i.e., $J^2 = J$).
\end{proposition}
\begin{proof}
We assume with no loss of generality that $S \subseteq A$ and
$\pi|_S = \id_S$.
Let $\sum_i x_i \otimes y_i$ be a separability element in $A \otimes_S A$.
Let $e = \sum_i \pi(x_i)y_i$.  Then $e$ satisfies $ex = \pi(x)e$ for all $x
\in A$.  Since $\pi(e) = 1$,
it follows that $e$ is idempotent.  Similarly, $f = \sum_i x_i \pi(y_i)$
satisfies $xf = f\pi(x)$, $\pi(f) = 1$ and is idempotent.  Then $$e =
\pi(f)e = ef = f.$$
Then $xe = e\pi(x)e = ex$ and $e$ is central.
Then $J = (1 - e)A$, whence $J$ is idempotent.
\end{proof}
An example of a ring epi splitting  in the next corollary would be the
one implicit in the
Wedderburn Principal Theorem for finite dimensional algebras.
\begin{corollary}
If $A/S$ is a split, separable extension with splitting map $\pi: A \to S$
a ring epimorphism with nilpotent kernel $J$,
then $J = 0$ and $A = S$.
\end{corollary}
Indeed a separable finitely generated (f.g.) extension of a separable
algebra is itself
separable \cite{HS}. The next proposition builds new separable extensions
from old
using
multiplicative bimodules \cite{P}.
\begin{proposition}
Suppose $R/S$ is a ring extension and $I$ is a multiplicative
$R$-bimodule.
Then $R/S$ is a separable extension
if and only if $A = R \oplus I$ is a separable extension of $T = S \oplus I$.
\end{proposition}
\begin{proof}
($\Rightarrow$)
Let $f = \sum_i x_i \otimes y_i \in R \otimes_S R$ be a separability element
for $R/S$.  Let $e$ be its image in $A \otimes_T A$ induced by $R
\hookrightarrow A$.
If $x \in I$, then:
$$ xe = \sum_i 1 \otimes xx_i y_i = \sum_i x_iy_i x \otimes 1 = ex. $$
We easily conclude that $e$ is a separability element for $A/T$.
($\Leftarrow$)  Since $I$ is an ideal in both $A$ and $T$,
the canonical epimorphism $\pi:  A \rightarrow R$ sends the separable
extension $A/T$ onto a separable extension $R/S$ \cite{HS}.
\end{proof}

As our final preliminary topic we recall Frobenius and QF
extensions.\footnote{\cite{K60,O,Mu64,Mu65,M65,M67,FMS,K}}
A ring extension $R/S$ is
\textit{Frobenius} if $R_S$ is  f.g. projective and $R \cong R^*$ as
$S$-$R$-bimodules:
note that this extends the notion of Frobenius algebra.
We recall also the Morita characterization
of Frobenius extensions \cite{M65}:
an extension $R/S$ is Frobenius iff induction and co-induction
(of $S$-modules to $R$-modules) are naturally isomorphic (cf.\ \cite{K}).

A ring extension $R/S$ is a left \textit{Quasi-Frobenius} (QF) extension
if ${}_SR$ is finitely generated projective
and ${}_RR_S$ is isomorphic to a direct summand of a finite direct product of
${}^*_RR_S$ with itself.  Equivalently, $R_S$ and ${}_SR$ are finitely
generated projectives
and ${}_SR^*_R$ is a direct summand of a finite direct sum of copies of
${}_SR_R$.  We similarly define right QF extensions \cite{Mu64,Mu65}.
There is no published example
of a right QF extension that is not left QF.
\end{section}

\begin{section}{Biseparable Bimodules}\selabel{2}
In this section we pose our question in the  more general terms of bimodules
rather than ring extensions.  There are two reasons for this.
First,  the problem has a more attractive symmetrical formulation
in terms of bimodules.  Second, Morita has shown in \cite{M67}
how to generate interesting examples of ring extensions from
bimodules via the endomorphism ring.

Let $R$ and $T$ be rings.  Given a bimodule ${}_T M_R$, there is a natural
$T$-bimodule
homomorphism,
$$\mu_M:\ M \otimes _R {}^*M \longrightarrow T, \ \ m \otimes f \mapsto
(m)f .$$
We next recall the definition of a separable bimodule  \cite{S71}.

\begin{definition}\delabel{2.1}
$M$ is separable, or $T$ is $M$-separable over $R$, if $\mu_M$ is a split
$T$-epimorphism.
\end{definition}

It follows trivially that $M_R$ is a generator module \cite{AF}.
By applying a splitting map to $1_T$,
we note that $M$ is separable iff there is an element $$e = \sum_i m_i
\otimes f_i \in M \otimes _R {}^*M,$$
called an \textit{$M$-separability element}, which satisfies $\mu_M(e) =
1_T$ and
$te = et$ for all $t \in T$.   As is the case with separability
elements and idempotents \cite{P,DMI}, $M$-separability elements
are usually not unique.

Retaining this notation, we recall a useful proposition and its proof
\cite[Proposition 3]{S71}.
But first a lemma which does not require $M$ to be separable.
\begin{lemma}
If $M_R$ is finitely generated projective, then $\alpha_M: M \otimes_R {}^*
M \to \Hom({}_RM^*, {}_R^*M)$
given by $$m \otimes f \longmapsto (g \mapsto g(m)f)$$
is a $T$-bimodule isomorphism.  Similarly, ${}_TM$ f.g. projective
implies that $$ \eta_M:  M^* \otimes_T M \to \Hom({}^*M_T, M^*_T), \ \ f
\otimes m \longmapsto
(h \mapsto f[(m)h])$$
is an $R$-bimodule isomorphism.
\end{lemma}

The proof of this and the next two propositions are left to the reader.
Now define $\beta_M:  \Hom({}_RM^*, {}_R^*M) \to T$ by
$$\beta_M(G) = \sum_i (x_i)[(f_i)G].$$

\begin{proposition}
\label{commutative triangle}
We have $\mu_M = \beta_M \circ \alpha_M$; i.e., the diagram
below is commutative. Whence
$M$ is separable iff there is an $R$-$T$-bimodule homomorphism $\gamma_M:
M^* \to {}^*M$
such that $(x_i)[(f_i)\gamma_M] = 1_T$.
\end{proposition}
$$\begin{diagram}
M \otimes_R {}^*M&&\rTo^{\cong}_{\alpha_M} && \Hom_R({}_RM^*, {}^*_RM)\\
&\SE_{\mu_M}&&\SW_{\beta_M}&\\
&&T&&
\end{diagram}$$

>From $\Hom_R(M^*, {}^*M)^T = \Hom_{R-T}(M^*, {}^*M)$ we note that
$\gamma_M$ corresponds under $\alpha_M$ to an $M$-separability element.

\begin{definition}
A bimodule ${}_T M_R$ is said to be \textit{biseparable} if $M$ and $M^*$
are separable
and ${}_TM$, $M_R$ are finite projective modules.
\end{definition}

We derive some consequences of assuming ${}_TM$ f.g. projective
and $M^*$ separable. First, since $M_R$ is reflexive, it follows that
${}^*(M^*) \cong M$
via the ``evaluation mapping'' from $M$ to $\Hom_R(M^*, R)$.  It follows that
$\mu_{M^*}: M^* \otimes_T M \to R$ under this identification
is the evaluation map given by $f \otimes m \mapsto f(m)$.
Let $\{ g_j \} \subset {}^*M$, $\{ y_j \} \subset M$
be  finite dual bases for the f.g. projective module ${}_TM$. We have the
following
easy analog of Proposition~\ref{commutative triangle}.

\begin{proposition}
The triangle below is commutative. Whence if $M^*$ is separable
there is a $R$-$T$-bimodule homomorphism $\rho_M: {}^*M \to M^*$
such that $\sum_j \rho_M(g_j)(y_j) =1_R$.
\end{proposition}
$$\begin{diagram}
M^* \otimes_T M & & \rTo{\cong}_{\eta_M} &  & \Hom({}^*M_T, M^*_T) \\
& \SE_{\mu_{M^*}} & & \SW & \\
& & R & &
\end{diagram}
$$
The downward map to the right is
given by $G \mapsto \sum_j G(g_j)(y_j)$.

For example, a bimodule ${}_TM_R$ yielding a Morita equivalence of $T$ and $R$
is biseparable, since ${}^*M \cong M^*$ as $R$-$T$-bimodules \cite{M67}
while $\mu_M$ and $\mu_{M^*}$ are isomorphisms.
There have been various studies of properties shared by rings $R$ and $T$
related by a bimodule ${}_TM_R$ in a Morita context and generalizations
of this \cite{AF,C,M,M65}.  A precursor of these studies
is the theorem of D.G. Higman \cite{Hi} that a finite group has finite
representation
type (f.r.t.) in characteristic $p$ iff its Sylow $p$-subgroup is cyclic,
which later became a corollary of the theorem of J.P. Jans \cite{Ja}
that for Artinian algebras $R \subseteq T$ in a split separable
extension, $R$ has f.r.t. iff $T$ has f.r.t. (cf.\ \cite{P}).
It is in this spirit
that the next theorem offers a  sample of shared properties of $R$ and $T$
linked by a biseparable bimodule ${}_TM_R$.

\begin{theorem}
If ${}_TM_R$ is biseparable, then
\begin{enumerate}
\item (Cf.\ \cite[Prop.\ 5, Theorem 2]{S71}.) $T$ is a QF ring if and only
if $R$ is a QF ring;
\item (Cf.\ \textit{Ibid.}) $T$ is semisimple if and only
if $R$ is semisimple;
%\item $T$ is left perfect if and only if $R$ is left perfect
\item weak global dimension $D(R) = D(T)$.
\end{enumerate}
\end{theorem}
\begin{proof}
(QF).
Assume $T$ is QF and $P_R$ is injective.  By the Faith-Walker theorem,
it will suffice to show that $P_R$ is projective.  Since
${}_T M$ is projective, then flat, we note that $H_T := \Hom_R(M,P)$
is injective, then projective.  Since $M_R$ is projective, we note
that $H \otimes_T M_R$ is projective.
But the evaluation mapping $$ ev:\ H \otimes_T M_R \to P_R,\ \ f \otimes m
\mapsto f(m)$$ is a split epi,
for if $\sum_j f_j \otimes m_j$ is an $M^*$-separability element,
then $p \mapsto pf_j \otimes m_j$ defines a splitting $R$-monic,
where $pf_j:\, m \mapsto pf_j(m)$.
Hence $P_R$ is isomorphic to a direct summand in $H \otimes_T M_R$
and is projective.

Assuming that $R$ is QF, and ${}_TQ$ is injective, we argue similarly
that ${}_RH' := \Hom_T(M,Q)$ is injective-projective and that
$Q$ is isomorphic to direct summand in the projective $T$-module
$M \otimes_R H'$.

(SEMISIMPLICITY). Suppose $T$ is semisimple and $P_R$ is a module.  It
suffices to note
that $P_R$ is projective. Since $H_T := \Hom_R(M,P)$ is projective and the
map $ev$ defined as above
is a split $R$-epimorphism, it follows that $P_R$ is isomorphic to a direct
summand
of the projective module $H \otimes_T M_R$.
Similarly, we argue that given $R$ semisimple and module ${}_TQ$,
${}_TQ$ is projective.

(WEAK GLOBAL DIMENSION.)  If $X_. \to N_T$ is a projective resolution of $N_T$,
then $X_. \otimes_T M_R \to N \otimes_T M_R$ is a projective resolution as well
since ${}_TM$ is flat and $M_R$ is projective.  Recall that
${\rm Tor}^T_n\, (N , Q)$
is the $n$'th homology group of the chain complex $X_. \otimes_T Q$ for each
non-negative integer $n$, so ${\rm Tor}^R_n\, (N \otimes_T M , {}^*M
\otimes_T Q)$
is the $n$'th homology group of $X_. \otimes_T M \otimes_R {}^* M \otimes_T
Q$.
At the level of chain complex, there is a split epi
$$ X_n \otimes_T M \otimes_R {}^* M \otimes_T Q \longmapsto X_n \otimes_TQ,\
\ x \otimes m \otimes g \otimes q \mapsto x(mg) \otimes q $$
which implies that ${\rm Tor}^T_n\, (N , Q)$ is isomorphic to a direct summand
in ${\rm Tor}^R_n\, (N \otimes_T M , {}^*M \otimes_T Q)$ for each $n$.
This shows that $D(R) \geq D(T)$.  A similar argument with left modules
shows that $D(T) \geq D(R)$.
\end{proof}

We remark that in trying to prove other shared homological properties
of biseparable ${}_TM_R$, particularly
one-sided notions, one may run into
the following complications: although the modules
${}^*M_T$ and ${}_R M^*$ are (quite easily seen to be)
f.g. projective, one should avoid assuming the
same of ${}^*_RM$ and $M^*_T$.

We next recall
the definition of Frobenius bimodule \cite{AF,K}.

\begin{definition}\delabel{3.8}
A bimodule ${}_TM_R$ is Frobenius if $M_R$, ${}_T M$ are f.g. projective
and ${}^*M \cong M^*$ as $R$-$T$-bimodules.
\end{definition}

Based on the many examples in \cite{K95} and
elsewhere, we propose the following problem, which
turns out to be almost equivalent
to  the ring extension formulation in the title:

\begin{problem}
\label{main}
Is a biseparable bimodule ${}_TM_R$ a Frobenius bimodule?
\end{problem}

For example, can we choose $\gamma_M$ and $\rho_M$ such
that they are inverses to one another?  The problem above
subsumes many interesting questions in various restricted
cases.  For example, what can be said
for the problem above if $T$ and $R$ are finite dimensional algebras?
There is an  affirmative answer in the next section
if one algebra is separable.

A generalization of Frobenius bimodule is a \textit{twisted Frobenius
bimodule}
${}_{\alpha}M_{\beta}$ where $\alpha: T \to T$
and $\beta: R \to R$ are ring automorphisms,
and the bimodule structure is now given by $t \cdot m \cdot r :=
\alpha(t)m\beta(r)$ (for the definition see \cite{K99}). We might
ask more widely
\begin{problem}
\label{twisted version}
Is a biseparable bimodule  a twisted
Frobenius bimodule?
\end{problem}

However, this problem is the same as the previous one if the following question
has an affirmative answer:

\begin{problem}
If a twisted Frobenius bimodule ${}_{\alpha}M_{\beta}$ is biseparable,
does this imply that  ${}_{\alpha}M_{\beta} \cong {}_TM_R$?
\end{problem}

We say that a twisted Frobenius bimodule is nontrivial if it not isomorphic
to an
untwisted Frobenius bimodule; for a $\beta$-Frobenius extension $R/S$
nontriviality means
that $\beta:S \to S$ is not an \textit{extended inner automorphism} in the
sense
that there is a unit $u \in R$ such that $\beta$ is conjugation by $u$
\cite{NT}.
We pose the last question since we have never observed a nontrivial
$\beta$-Frobenius extension
(e.g. in \cite{FMS, NT}) which was
simultaneously split and separable (cf.\ next section).
In this more limited setting, which covers Hopf subalgebras of finite
dimensional
Hopf algebras, the question becomes:

\begin{problem}
If a $\beta$-Frobenius extension is split and separable, is $\beta$
an extended inner automorphism?
\end{problem}
We will return to a discussion of this problem in the next section.

\end{section}

\begin{section}{Biseparable Extensions}\selabel{3}
Suppose $R/S$ is a ring extension.  Letting $M = {}_R R_S$ in the definition
of separable bimodule, we observe the following lemma \cite{S71}.
\begin{lemma}
$R/S$ is a separable extension iff ${}_R R_S$ is a separable bimodule.
\end{lemma}
Dually, we let $M = {}^*({}_R R_S) \cong {}_SR_R$ and observe the following
lemma \cite{K96}.
\begin{lemma}
$R/S$ is a split extension iff ${}_SR_R$ is a separable bimodule
iff $({}_RR_S)^*$ is a separable bimodule.
\end{lemma}
The $M$-separability element in this case is a \textit{bimodule projection}
$E:  {}_SR_S \to {}_SS_S$, which implies $R/S$ is a proper extension.
$E$ is also called a \textit{conditional expectation} if it satisfies
additional properties in subfactor theory.
>From the last two lemmas, it follows that:
\begin{lemma}
\label{lemma-def}
$R/S$ is a split, separable, two-sided finite projective extension iff
${}_R R_S$ and ${}_S R_R$ are biseparable bimodules iff ${}_RR_S$
is biseparable and ${}_SR$ is f.g. projective iff
${}_SR_R$ is biseparable and $R_S$ is f.g. projective.
\end{lemma}
We call $R/S$ a \textit{biseparable extension} if
any of the equivalent conditions
in Lemma~\ref{lemma-def} are satisfied.
Additionally, we have the following lemma \cite{K}.
\begin{lemma}
$R/S$ is a Frobenius extension iff  ${}_RR_S$ or  ${}_SR_R$ is
a Frobenius bimodule.
\end{lemma}
The last two lemmas lead to the seemingly restricted formulation of
Problem~\ref{main},
also the title of this article.
\begin{problem}
\label{title}
Are biseparable extensions Frobenius?
\end{problem}

Surprisingly, this problem is almost equivalent to Problem~\ref{main} because
of the endomorphism ring theorems for Frobenius bimodules \cite[Chapter 2]{K}.
Suppose we knew an affirmative answer to the somewhat weaker problem
where biseparable extension includes \textit{left} f.g. projective,
split separable extensions.
Given a biseparable bimodule ${}_TM_R$, we know from Sugano
that $\mathcal{E} = \End(M_R)$ is a left f.g. projective, split, separable
extension
of $T$ (whose elements are identified in $\mathcal{E}$ with left multiplication
operators) \cite[Theorem 1, Prop.\ 2]{S71} (cf.\ \cite[Theorem 3.1]{K99}).
Then $\mathcal{E}/T$ is a Frobenius extension by our affirmative answer to
the weak Problem~\ref{title}.  Since ${}_{\mathcal{E}}M_R$ is faithfully
balanced by Morita's Lemma, it follows from \cite[Theorem 1.1]{M67}
that ${}_{\mathcal{E}}^*M_R \cong M^*$ and then from
the endomorphism ring theorem-converse \cite[Theorem 2.8]{K}
that ${}_TM_R$ is a Frobenius bimodule.

What evidence do we have then for proposing Problem~\ref{title}?
First, if $R/S$ is an $S$-algebra, we are in the situation of
a faithfully projective separable algebra, which is Frobenius by the
Endo-Watanabe
Theorem \cite{EW} discussed in the introduction.
This implies by elementary considerations that $k$-algebra extensions
of the form $R \otimes_k A$ over $R$ are Frobenius
if $A$ is faithfully projective $k$-separable.

Second, there are the many examples of split, separable, Frobenius
extensions \cite{K95,K96} and apparently none that contradict in the
literature for
noncommutative rings and ring
extensions \cite{JS,L,M67,NT}. We recall from \cite{K95,K96} that some of
the examples of
$R/S$ split, separable Frobenius are the following:
\begin{enumerate}
\item Let $R = E$ a field of characteristic $p$
(zero or prime) and $S  = F$ a subfield such that $E/F$ is
a finite separable extension where $p$ does not divide $[E:F]$.
The classical trace map from $E$ into $F$ is a Frobenius homomorphism
in this example.

\item  Let $R$ be the group algebra $k[G]$ for  a discrete group $G$, $k$ a
field, and $S = k[H]$,
where $H$ is a subgroup of $G$ with finite index not divisible by the
characteristic of $k$ (e.g., char $k = p$ and $H$
a Sylow $p$-subgroup).  Note that if the characteristic of $k$
divides $[G:H]$ we have an example of a split Frobenius extension
which is \textit{not} separable.  By \cite{S71}
the endomorphism ring extension $\End(R_S)/R$ is a separable
Frobenius extension which is \textit{not} split.

\item Let $R$ and $S$ be  algebras over a commutative ring $k$ such
that $R/S$ is an Hopf-Galois extension over the f.g. projective Hopf
$k$-algebra $H$ which is separable and coseparable over $k$.

\item Let $R$ be a type $II_1$ factor, $S$ a subfactor of $R$ of finite
Jones index, as discussed in the introduction.

\end{enumerate}

Third,  Sugano's result \cite{S70} for when H-separable extensions
are Frobenius  is evidence for biseparable implies
Frobenius. This is because an H-separable extension is a strong
type of separable extension \cite{H}:  see Section~4 for a separable
extension which is not H-separable. Thus the
result that a (one-sided) f.g. projective split H-separable extension
is a (symmetric) Frobenius extension is a particular case
of a biseparable extension which is Frobenius
(cf.\ \cite[Section 2.6]{K}\footnote{\cite[Theorem2.25]{K}
should read ``Suppose $A/S$ is a right progenerator split H-separable...''
as the proof clearly shows.}).

The next proposition shows that a biseparable extension is almost
a two-sided QF extension in a certain sense.
If a module $M_R$ is isomorphic to a direct summand in another module
$P_R$, we denote this by $M_R <_{\oplus}\ P_R$.
\begin{proposition}
Suppose $R/S$ is a biseparable extension.  Then all $R$-modules are
$S$-relative injective and $S$-relative projective; moreover,
$ R^*_R $ and ${}^*_R R$ are generator modules.
\end{proposition}
\begin{proof}
The first statement follows from the fact that a separable
extension is both right and left semisimple extension
and properties of these \cite{HS}.

For the second statement, we first establish an interesting isomorphism
below involving $R$ and its dual $R^*$.
On the one hand, since $R/S$ is a separable extension, $\mathcal{E} =
\Hom(R_S,R_S)$ is a split
extension of $R$, for if $\sum_i x_i \otimes y_i$ is a separability
element we define a bimodule projection by
$\mathcal{E} \to R$ by $f \mapsto \sum_i f(x_i)y_i$ (cf.\ \cite{Mu65}).
Then as $R$-bimodules,
$\mathcal{E} \cong R \oplus M$ for some  $M$:
moreover, by restriction this is true as $S$-$R$-bimodules.  On the other
hand, since $R/S$ is split,
it follows that for some $S$-bimodule $N$,
which is left and right projective $S$-module,
$R \cong S \oplus N$ as $S$-bimodules; whence
$$\mathcal{E}
\cong R^* \oplus \Hom(R_S,N_S)$$
as $S$-$R$-bimodules.
Putting together the two isomorphisms for $\mathcal{E}$, we obtain
\begin{equation}
\label{isomorphism}
R \oplus M \cong \Hom(R_S,S_S) \oplus \Hom(R_S,N_S).
\end{equation}

Since $N_S$ is f.g. projective, there is a module $P_S$ such that
$N_S \oplus P_S \cong S^r_S$.  Then applying $\Hom_S(R_S,-)$ to this:
$$
\Hom(R_S,N_S) \oplus \Hom(R_S,P_S) \cong \Hom(R_S, S_S)^r.
$$
Combining this with Eq. (\ref{isomorphism}), we obtain
\begin{equation}
\label{iso2}
R_R \oplus M_R \oplus \Hom(R_S,P_S) \cong R^{* \, r+1}_R.
\end{equation}
This establishes that $ R_R <_{\oplus}\ R^{*\, r+1}_R$.

We similarly conclude ${}_R R <_{\oplus} {}^*_R R^{t+1}$ by combining
the split extension $\mathcal{E}' := \Hom({}_S R, {}_S R)/ R$ with
the $S$-bimodule
isomorphism $R \cong S \oplus N$ and the existence of ${}_S Q$
such that ${}_SQ \oplus {}_SN \cong {}_S S^t$.
\end{proof}
>From the proof just completed,
we obtain a corollary worth noting for its
relatively easy proof. ${}_AM_A$ is said to be {\it centrally
projective}\index{bimodule!centrally projective}
if $ {}_AM_A <_{\oplus}\ {}_AA^n_A$ for some positive integer $n$.

\begin{corollary}(Cf.\ \cite[Theorem 2]{S70})
If $R/S$ is centrally projective biseparable extension, then it is a QF
extension.
\end{corollary}
\begin{proof}
Since there is a bimodule ${}_SP'_S$ such that ${}_S R_S \oplus {}_SP'_S
\cong {}_S S^r_S$,
we combine this with ${}_S R_S \cong {}_S S_S \oplus {}_S N_S$ to see that Eq.\
(\ref{iso2}) is an $S$-$R$-isomorphism, whence $R/S$ is a right QF extension
\cite{Mu64}.  Similarly, we show $R/S$ to be a left QF extension.
\end{proof}

The proposition and corollary\footnote{In terms of finite depth
for subfactors \cite{JS} and ring extensions \cite{NK, SzK}, the corollary
states that a depth one biseparable extension is QF.  A
 depth two biseparable extension is QF as well \cite{SzK},
where depth two is the condition that $R \otimes_S R \oplus * \cong \oplus^n R$
as $S$-$R$ and $R$-$S$-bimodules.} above lead naturally to the problem below,
a weakening of Problem~\ref{title}.
\begin{problem}
\label{weak version}
Are biseparable extensions QF?
\end{problem}

The next theorem provides a solution of Problem~\ref{title}
in case
$A$ or $S$ is a separable algebra.  We assume our algebras
to be faithful.

\begin{theorem}
Suppose $A/S$ is a biseparable extension of $k$-algebras
with $k$ a commutative ring.
If either  $S$ is a $k$-projective separable $k$-algebra, or
$A$ is a separable $k$-algebra with $k$ a field,
then $A/S$ is a Frobenius extension.
\end{theorem}
\begin{proof}
The proof does not make use of ${}_SA$ being f.g. projective.
Suppose $k$ is a field  and $A$ is $k$-separable.
Let $D(R)$ denote the right global dimension of a ring $R$
and $d(M)$ denote the projective dimension of a module $M_R$. Then $D(A) = 0$
since $A$ is finite dimensional semisimple.  By Cohen's Theorem
for split extensions \cite{Kap},
$$ D(S) \leq D(A) + d(A_S) $$
whence $D(S) = 0$ and $S$ is semisimple \cite{Kap}.  Then $A$ and $S$
are finite dimensional semisimple algebras.
It follows from \cite{EN} that  $S$ and $A$ are symmetric algebras.

Similarly, we arrive at symmetric algebras $S$ and $A$ via \cite{EW} under
the assumption
that $S$ is $k$-projectively $k$-separable with no restriction on $k$.
For then $A$ is $k$-projective
and $k$-separable by transitivity for projectivity and separability.

Now we compute using the
\textit{bimodule isomorphisms} $A \cong A^*$ and $S \cong S^*$
and the  hom-tensor adjunction:
\begin{eqnarray*}
{}_SA_A & \cong & {}_S \Hom_k( A,k)_A \\
& \cong & \Hom_k(A \otimes_S S_{S}, k)_A \\
& \cong & {}_S\Hom_S(A_S, \Hom_k(S,k)_S)_A \\
& \cong & {}_S \Hom_S(A_S,S_S)_A.
\end{eqnarray*}
Then, since $A_S$ is f.g. projective, $A/S$ is a Frobenius extension.
\end{proof}

Part of the theorem is true without the hypothesis
of biseparable extension for a finite-dimensional \textit{Hopf subalgebra}
pair $H \supseteq K$:
if $H$ is semisimple, then $K$ is semisimple \cite[2.2.2]{M},
and $H/K$ is a Frobenius extension (cf.\ \cite[1.8]{FMS}).

Finally, Problem (\ref{title}) can be widened to twisted extensions,
as Problem (\ref{main}) was widened to twisted Frobenius bimodules
in Problem (\ref{twisted version}).

\begin{problem}
Are biseparable extensions $\alpha$-$\beta$-Frobenius?
\end{problem}
We refer the reader to \cite{M65,K99} for the definition of these
twisted extensions, which are more general than the usual $\beta$-Frobenius
extensions.  Also, Problem (\ref{weak version}) has a twisted enlargement.
%The authors do not know of a f.g. free separable extension that is
%not $\beta$-Frobenius, and this would naturally lead to another question.

\end{section}
\begin{section}{Examples and Counterexamples}\selabel{4}
In this section, we consider weakening the definition of biseparability
in various ways, and find examples of non-Frobenius extension for
each such case. We will see an example of non-finitely generated
projective separable extension, which  is  an obstruction
to extending Villamayor's theorem \cite[Prop.\ 10.3]{P} and Tominaga's
theorem \cite{T}.

\begin{lemma}
Suppose $k$ is a commutative ring
and $R$ is a $k$-algebra with $xy = 1$ but $yx \not= 1$.
Then $R$ is a separable
extension over $S = k1 + yRx$.
\end{lemma}
\begin{proof}
We note that $e = yx$ is a nontrivial idempotent in $S$.
Consider $f = x \otimes y \in R \otimes_S R$. Of course, $\mu(f) = 1$.
We compute with $r \in R$:
$$ rx \otimes y = xyrx \otimes y = x \otimes yrxy = x \otimes yr. \qed $$
\renewcommand{\qed}{}\end{proof}

Now if $R$ is a \textit{finitely} generated, projective
$k$-algebra, it is well-known that $xy = 1$ implies $yx = 1$.
So we let $V$ be a countably infinite rank free
$k$-module and
$$R = \End_k(V) = \{  \left( \begin{array}{cc}
a & \underline{b} \\
\underline{c}^t & D
\end{array}
\right) |\, a \in k, D \in M_{\infty}(k), \underline{b}, \underline{c}^t
\in M_{1 \times \infty  }(k) \}, $$ where $\underline{c}^t$ and $D$ are
column-finite and
$X \mapsto X^t$ denotes transpose.
$R$
is a ring satisfying the hypothesis in the lemma with elements $x, y \in R$
given in terms of the matrix units $e_{i,j}$ by
$$ x = \sum_{n=1}^{\infty} e_{n,n+1}, \ \ \ y = \sum_{n=1}^{\infty}
e_{n+1,n}. $$
Clearly, $xy = 1$ but $$e = yx = \left( \begin{array}{cc}
0 & \underline{0} \\
\underline{0}^t &  I
\end{array}
\right) $$
and the $k$-subalgebra,
$$ S = k1_R + yRx = \{ \left( \begin{array}{cc}
a & \underline{0} \\
\underline{0}^t & D
\end{array}
\right) |\, a \in k, D \in M_{\infty}(k) \} $$

\begin{proposition}
$R/S$ is a split, separable,  projective and non-finitely generated
extension.
\end{proposition}
\begin{proof}
We have seen in the lemma that $R/S$ is separable.  It is split since
we easily check that $E: R \to S$ below is a bimodule projection:
$$
E \left( \begin{array}{cc}
a & \underline{b} \\
\underline{c}^t & D
\end{array}
\right) =  \left( \begin{array}{cc}
a & \underline{0} \\
\underline{0}^t & D
\end{array}
\right)  $$

$R_S$ is countably generated projective, since $f_n: R_S \to S_S$ defined by
$$f_1 \left( \begin{array}{cc}
a & \underline{b} \\
\underline{c}^t & D
\end{array}
\right) =  \left( \begin{array}{cc}
a & \underline{0} \\
\underline{0}^t & \left( \begin{array}{c} \underline{b} \\ 0 \end{array}
\right)
\end{array}
\right)\ \  f_2 \left( \begin{array}{cc}
a & \underline{b} \\
\underline{c}^t & D
\end{array}
\right) =  \left( \begin{array}{cc}
0 & \underline{0} \\
\underline{0}^t & D
\end{array}
\right) $$
and $ f_{2+n}\left( \begin{array}{cc}
a & \underline{b} \\
\underline{c}^t & D
\end{array}
\right) = \left( \begin{array}{cc}
c_n & \underline{0} \\
\underline{0}^t & 0
\end{array}
\right) $, where $\underline{c} = (c_1, c_2,\ldots)$,
which satisfies the dual bases equation,
$$ A =  \left( \begin{array}{cc}
a & \underline{b} \\
\underline{c}^t & D
\end{array}
\right)  =  \left( \begin{array}{cc}
1 & (1\ \underline{0}) \\
\underline{0}^t & 0
\end{array}
\right) f_1 (A) +  \left( \begin{array}{cc}
0 & \underline{0} \\
\underline{0}^t & I
\end{array}
\right)f_2 (A) + \sum_{n=1}^{\infty} e_{n+1,1}f_{2+n} (A)  $$
By tranposing the elements and maps above, we similarly find a countable
projective
base for ${}_SR$.  The rest of the proof is now clear.
\end{proof}

\begin{corollary}
$R/S$ is not  Frobenius  and not H-separable.
\end{corollary}
\begin{proof}
For $R/S$ to be a Frobenius extension,
we must have $R_S$  finitely generated from the very start.  A
  right
projective H-separable extension is right f.g. by \cite{T}.
\end{proof}

Finally, we consider various weakenings of Problem~\ref{title}
and note that they all have known counterexamples.
There are easy examples of  separable  extensions which
are not Frobenius, such as the rationals $\mathcal{Q}$ extending
the integers $\mathcal{Z}$.  There are even f.g. free separable
extensions that are not Frobenius in the ordinary sense, but
are $\beta$-Frobenius \cite{K99}. As for a f.g. split, separable extension
that is not  Frobenius, here is one that is not
flat, whence not projective:  consider
$\mathcal{Z} \oplus \mathcal{Z}_2$ as a $\mathcal{Z}$-algebra,
which has projection onto its first factor as splitting bimodule projection,
is a direct sum of separable algebras -- whence separable -- and its second
factor
is of course not flat over $\mathcal{Z}$.

The following example of a f.g. projective H-separable (therefore separable
\cite{H1})
extension $A/S$ which is not even twisted QF is certainly worth a mention:
let $A$ be the $n$ by $n$ matrices over a field with $S$ the upper
triangular matrices (including diagonal matrices).  Although this
example is not Frobenius in any sense that we have mentioned in this paper
(left as an exercise to the reader), we note that
$\End A_S \cong A \otimes_S A (\cong A)$. However, there is no
published example of a one-sided progenerator H-separable extension
which is not Frobenius.

There are clearly many examples of split
extensions that are not Frobenius, let alone f.g.
Asking for a split, f.g. projective extension that fails to be
Frobenius is not hard:  for example, let $R$ be the upper triangular
$n \times n$ matrix algebra with splitting $S \oplus I$ where
$S$ is the subalgebra of diagonal matrices and $I^n = 0$.
It is well-known that $R$ is not a QF-algebra, and certainly not Frobenius
\cite{L},
but $S$ is semisimple, so $R_S$ and ${}_SR$ are f.g. projective; moreover,
$R/S$ cannot be Frobenius by the transitivity property
of Frobenius extension (cf.\ \cite{SK2}). This is an example too of
$\Hom(R_S,S_S)$ not being a projective right $R$-module.

As a last cautionary example, we consider $R = \mathcal{Z}_2 \oplus
\mathcal{Z}_2$
and $S = \mathcal{Z}_2$.  It is easy to check that $R/S$ is split,
separable and
Frobenius, even f.g. free.  But there are only two bimodule projections
$E: R \to S$, neither of which is a Frobenius homomorphism, i.e. in possession
of dual bases \cite{K}.  The Frobenius homomorphism in this example is
unique, since
the group of units in $C_R(S)$ consists only of the identity \cite{K}.
\end{section}

\begin{section}{Categorical interpretation}\selabel{5}
Let $F:\ {\mathcal C}\to {\mathcal D}$ be a contravariant functor. $F$ induces
a natural transformation ${\mathcal F}:\ \Hom_{\mathcal C}(\bullet,\bullet)\to
\Hom_{\mathcal D}(F(\bullet),F(\bullet)),~~{\mathcal F}_{C,C'}(f)=F(f)$.
$F$ is called a separable functor \cite{NVO} if ${\mathcal F}$ splits,
i.e. there exists a natural transformation ${\mathcal P}:\
\Hom_{\mathcal D}(F(\bullet),F(\bullet))\to \Hom_{\mathcal C}(\bullet,\bullet)$
such that ${\mathcal P}\circ {\mathcal F}$ is the identity natural
transformation
on $\Hom_{\mathcal C}(\bullet,\bullet)$.

\begin{proposition} \cite{RAFAEL}
Assume that $F$ has a right adjoint $G$, and let $\eta:\ 1_{\mathcal C}\to
GF$ and $\varepsilon:\ FG\to 1_{\mathcal D}$ be the unit and counit of
the adjunction.\\
$F$ is separable if and only if there exists a natural transformation
$\nu:\ GF\to 1_{\mathcal C}$ such that $\nu\circ\eta$ is the identity natural
transformation on ${\mathcal C}$.\\
$G$ is separable if and only if there exists a natural transformation
$\zeta:\ 1_{\mathcal D}\to FG$ such that $\varepsilon\circ\zeta$ is the
identity natural
transformation on ${\mathcal D}$.
\end{proposition}

The terminology stems from the fact that, for a ring homomorphism
$i:\ R\to S$, the restriction of scalars functor is separable
if and only if $S/R$ is separable (see \cite{NVO}, \cite{RAFAEL});
Separable functors satisfy the following version of Maschke's
Theorem: if a morphism $f$ in ${\mathcal C}$ is such that $F(f)$
has a left or right inverse in ${\mathcal D}$, then $f$ has
a left or right inverse in ${\mathcal C}$. Separable functors have
been studied in several particular cases recently, see e.g.
\cite{CastanoGN97}, \cite{CastanoGN98}, \cite{CaenepeelIMZ99}.\\
The functor $F$ is called Frobenius if $F$ has a right adjoint
$G$ that is at the same time a right adjoint. We will then
say that $(F,G)$ is a Frobenius pair. Now the terminology is
inspired by the property that a ring homomorphism $i:\ R\to S$
is Frobenius if and only if the restriction of scalars functor
is Frobenius. Frobenius pairs were introduced in \cite{M65},
and studied more recently in \cite{CMZ97} and \cite{CGN99}.
For more details and examples of separable functors
and Frobenius functors, we refer the reader to \cite{CMZ}.\\
Suppose we know that $(F,G)$ is a Frobenius pair. Then Rafael's
Theorem can be simplified: we can give an easier criterion
for $F$ or $G$ to be separable. First we need a result on
adjoint functors. Let $\varepsilon$ and $\eta$ be
the counit and unit of an adjunction $(F,G)$, and recall that
\begin{equation}\eqlabel{5.1.1}
\varepsilon_{F(C)}\circ F(\eta_C)=I_{F(C)}~~{\rm and}~~
G(\varepsilon_D)\circ \eta_{G(D)}=I_{G(D)}
\end{equation}
for all $C\in {\mathcal C}$ and $D\in {\mathcal D}$.

\begin{lemma}\lelabel{5.1}
Let $(F,G)$ be an adjoint pair functors, then we have isomorphisms
$$\Nat(F,F)\cong\Nat(G,G)\cong\Nat(1_{\mathcal
C},GF)\cong\Nat(FG,1_{\mathcal D})$$
\end{lemma}

\begin{proof}
In fact this is an easy consequence of \cite[Section 1.15, Lemma 3 and
Section 2.1, Corollary 1 and Lemma]{Pareigis}. We restrict to giving
a description of the isomorphism $\Nat(G,G)\cong\Nat(1_{\mathcal C},GF)$,
the other isomorphisms are similar.
For a natural transformation $\theta:\ 1_{\mathcal C}\to GF$,
we define $\alpha=X(\theta):\ G\to G$ by
\begin{equation}\eqlabel{5.1.2}
\alpha_D=G(\varepsilon_D)\circ\theta_{G(D)}
\end{equation}
Conversely, for $\alpha:\ G\to G$, $\theta=X^{-1}(\alpha):\
1_{\mathcal C}\to GF$ is defined by
\begin{equation}\eqlabel{5.1.3}
\theta_{C}=\alpha_{F(C)}\circ \eta_C
\end{equation}
The fact that $X$ and $X^{-1}$ are each others
inverses is proved using \equref{5.1.1}.\end{proof}

\begin{corollary}\colabel{5.2}
Let $(F,G)$ be a Frobenius pair of functors, then we have isomorphisms
\begin{eqnarray*}
\Nat(F,F)\cong\Nat(G,G)&\cong&\Nat(1_{\mathcal
C},GF)\cong\Nat(FG,1_{\mathcal D})\\
&\cong&\Nat(GF,1_{\mathcal C})\cong\Nat(1_{\mathcal D},FG)
\end{eqnarray*}
\end{corollary}

For a Frobenius pair $(F,G)$, we will write $\nu:\ GF\to 1_{\mathcal C}$
and $\zeta:\ 1_{\mathcal D}\to FG$ for the counit and unit of the adjunction
$(G,F)$. For all $C\in {\mathcal C}$ and $D\in {\mathcal D}$, we then have
\begin{equation}\eqlabel{5.2.1}
F(\nu_C)\circ\zeta_{F(C)}=I_{F(C)}~~{\rm and}~~
\nu_{G(D)}\circ G(\zeta_D)=I_{G(D)}
\end{equation}

\begin{proposition}\prlabel{5.3}
Let $(F,G)$ be a Frobenius pair, and let $\eta,~\varepsilon,~\nu$ and $\zeta$
be as above. The following statements are equivalent:
\begin{itemize}
\item $F$ is separable;
\item $\exists\alpha\in \Nat(F,F):~\nu_C\circ G(\alpha_C)\circ\eta_C=I_C$
for all $C\in {\mathcal C}$;
\item $\exists\beta\in \Nat(G,G):~\nu_C\circ \beta_{F(C)}\circ\eta_C=I_C$
for all $C\in {\mathcal C}$.
\end{itemize}
We have a similar characterization for the separability of $G$: the
statements
\begin{itemize}
\item $G$ is separable;
\item $\exists\alpha\in \Nat(F,F):~\varepsilon_D\circ\alpha_{G(D)}\circ
\zeta_D=I_D$
for all $D\in {\mathcal D}$;
\item $\exists\beta\in \Nat(G,G):~\varepsilon_D\circ F(\beta_D)\circ
\zeta_D=I_D$
for all $D\in {\mathcal D}$.
\end{itemize}
are equivalent.
\end{proposition}

\begin{proof}
Assume that $F$ is separable. By Rafael's Theorem, there exists $\tilde{\nu}\in
\Nat(GF,1_{\mathcal C})$ such that $\tilde{\nu}_C\circ \eta_C=I_C$ for all
$C\in {\mathcal C}$. Let $\alpha:\ F\to F$ be the corresponding natural
transformation
of \coref{5.2}, i.e. $\alpha_C= F(\tilde{\nu}_C)\circ \zeta_{F(C)}$, and
$\tilde{\nu}_C=\nu_C\circ G(\alpha_C)$, and the first implication of the
Proposition follows. The converse follows trivially from Rafael's Theorem.
All the other equivalences can be proved in a similar way.
\end{proof}

\subsection*{Application to bimodules}
We use the notation of \seref{2}: let $R$ and $T$ be rings, and
$M$ a $(T,R)$-bimodule. We have already seen that $M^*=\Hom(M_R,R_R)$
and ${}^*M=\Hom({}_TM,{}_TT)$ are $(R,T)$-bimodules. As given
in the preliminaries, 
$\Hom({}_TM,{}_TM)$  and $\Hom(M_R,M_R)$
are respectively the natural $(R,R)$-bimodule and  $(T,T)$-bimodule. 
Furthermore, the induction functor
$$F=M\ot_R\bullet:\ {}_R{\mathcal M}\to {}_T{\mathcal M}$$
has a right adjoint
$$G=\Hom({}_TM,{}_T\bullet):\ {}_T{\mathcal M}\to {}_R{\mathcal M}$$
For $Q\in {}_T{\mathcal M}$, $\Hom({}_TM,{}_TQ)\in {}_R{\mathcal M}$ via
$(m)(rf)=(mr)(f)$. For $\kappa:\ Q\to Q'$ in ${}_T{\mathcal M}$,
we put $G(\kappa)(f)=\kappa\circ f$.
We also have a functor
$$G'={}^*M\ot_T\bullet:\ {}_T{\mathcal M}\to {}_R{\mathcal M}$$
and a natural transformation $\gamma:\ G'\to G$ given by
$$\gamma_Q:\ {}^*M\ot_T Q\to \Hom({}_TM,{}_TQ),~~
(m)(\gamma_Q(f\ot_T q))=(m)fq$$
We will also write $\gamma_Q(f\ot_T q)=f\cdot q$. $\gamma_Q$ is
well-defined on the tensor product, and left $R$-linear, so
\begin{equation}\eqlabel{5.4.1}
(rft)\cdot q=r(f\cdot (tq))
\end{equation}
for all $r\in R$ and $t\in T$. If $Q$ is a $(T,R)$-bimodule
(for example, $T=M$), then $\gamma_Q$ is also right $R$-linear, and
we have
\begin{equation}\eqlabel{5.4.2}
f\cdot (qr)=(f\cdot q)r
\end{equation}
Now assume that ${}_TM$ is finitely generated and projective,
and consider a dual basis $\{n_j\in M,~g_j\in {}^*M\}$, i.e.
$$m=\sum_j (m)g_jn_j$$
for all $m\in M$, or $I_M=\sum_j\gamma_M(g_j\ot_T n_j)$. Then
$\gamma:\ G\to G'$ is a natural isomorphism, and for all
left $T$-linear $f:\ M\to Q$, we have
$$\gamma_Q^{-1}(f)=\sum_j g_j\ot_T f(n_j)$$
In order to decide whether $F$ or $G$ is separable, or whether $(F,G)$
is a Frobenius pair, we have to investigate natural transformations
$GF\to 1_{{}_R{\mathcal M}}$ and $1_{{}_T{\mathcal M}}\to FG$. This is
done in the next two Propositions.
We thank the referee for pointing out to us that Proposition \ref{pr:5.5}
and \ref{pr:5.6} can be derived from \cite[Theorem 3.5.6]{Popescu}.

\begin{proposition}\prlabel{5.5}
Let $F$ and $G$ be as above, and consider
$$V=\Nat(GF,1_{{}_R{\mathcal M}}),~~V_1=\Hom({}_R\Hom({}_TM,{}_TM)_R,{}_RR_R)$$
$$V_2=\Hom({}_R({}^*M)_T,{}_R(M^*)_T)$$
Then we have maps
$$V\rTo^{\alpha} V_1\rTo^{\alpha_1} V_2$$
which are isomorphisms if $M$ is finitely generated and projective as a
left $T$-module.
\end{proposition}

\begin{proof}
For $\nu\in V$, we put $\alpha(\nu)=\nu_R$. By definition, $\nu_R$ is
left $R$-linear. For any $s\in R$, we consider the left $R$-linear map
$m_s:\ R\to R$, $m_s(r)=rs$. The naturality of $\nu$ implies that
$\nu_R(f)s=\nu_R(GF(m_s)(f))=\nu_R(fs)$, so $\nu_R$ is also
right $R$-linear.\\
$\alpha_1:\ V_1\to V_2$ is given by $\alpha_1(\ol{\nu})=\ol{\phi}$,
with
$$\ol{\phi}(f)(m)=\ol{\nu}(\gamma_M(f\ot m))=\ol{\nu}(f\cdot m)$$
Using \equref{5.4.1} and \equref{5.4.2}, we easily deduce that
$\ol{\phi}(f)$ is right $R$-linear, and that $\ol{\phi}$ is left
$R$-linear and right $T$-linear.\\
Assume that ${}_TM$ is finitely generated projective, and, as above,
assume that $\{n_j\in M,~g_j\in {}^*M\}$ is a dual basis. We can
then define the inverse $\alpha^{-1}$ of $\alpha$ as follows.
We view $\ol{\nu}\in V_1$ as a map $\ol{\nu}:\ {}^*M\ot_T M\to R$,
and we identify $G$ and $G'$. We then define $\nu\in V$ by
$$\nu_P=\ol{\nu}\ot_R P:\ GF(P)\cong {}^*M\ot_T M\ot_RP\to R\ot_RP\cong P$$
It is clear that $\nu$ is natural and that $\alpha$ and $\alpha^{-1}$
are each others inverses.\\
For $\ol{\phi}\in V_2$, we define $\alpha_1^{-1}(\ol{\phi})=\ol{\nu}$ by
$$\ol{\nu}(\varphi)=\sum_j \ol{\phi}(g_j)((n_j)\varphi)$$
for all $\varphi:\ M\to M$ in ${}_T{\mathcal M}$. Straightforward computations
yield that $\alpha_1^{-1}$ is well-defined, and is indeed an inverse
of $\alpha_1$.
\end{proof}

In a similar fashion, we have:

\begin{proposition}\prlabel{5.6}
Let $F$ and $G$ be as above, and consider
$$W=\Nat(1_{{}_T{\mathcal M}}, FG),~~W_1=\{e\in M\ot_R{}^*M~|~te=et,~{\rm
for~all~}
t\in T\}$$
$$W_2=\Hom({}_R(M^*)_T,{}_R({}^*M)_T)$$
We have maps
$$W\rTo^{\beta}W_1\rTo^{\beta_1}W_2$$
$\beta$ is an isomorphism, and $\beta_1$ is an isomorphism if $M$ is
finitely generated as a right $R$-module.
\end{proposition}

\begin{proof}
The proof is similar to the previous one, so we restrict to giving
the connecting maps.
For a natural transformation $\zeta:\ 1_{{}_T{\mathcal M}}\to FG$, we define
$e=\beta(\zeta)=\zeta_T(1)\in FG(T)=M\ot_R{}^*M$.
Conversely, take $e=\sum_i m_i\ot_R f_i\in W_1$. $\zeta\in W$
is defined as follows: for all $Q\in {}_T{\mathcal M}$, we put
$$\zeta_Q:\ Q\to M\ot_R \Hom({}_TM,{}_T,Q),~~
\zeta_Q(q)=\sum_i m_i\ot_R f_i\cdot q$$
Now we define $\beta_1:\ W_1\to W_2$. For $e=\sum_i m_i\ot_R f_i\in W_1$,
we let
$$\beta_1(e)=\phi:\ M^*\to {}^*M$$
be given by
$$\phi(h)=\sum_i h(m_i)f_i~~{\rm or}~~
(m)\phi(h)=\sum_i (mh(m_i))f_i$$
Straightforward computations show that $\phi(h)\in {}^*M$, and that
$\phi$ is $(R,T)$-bilinear.
\end{proof}

The two previous results can be used to decide when the induction functor
$F$ and its adjoint $G$ are separable or Frobenius. Let us first look
at separability.

\begin{corollary}\colabel{5.7}
Let $M$ be a $(T,R)$-bimodule, and assume that $M$ is finitely generated
projective as a left $T$-module, with finite dual basis $\{n_j,g_j\}$.
Then the following assertions are equivalent:
\begin{itemize}
\item $F=M\ot_R\bullet$ is a separable functor;
\item there exists $\ol{\nu}\in V_1$ such that $\ol{\nu}(I_M)=1_R$;
\item there exists $\ol{\phi}\in V_2$ such that
$\sum_j\ol{\phi}(g_j)(n_j)=1_R$.
\end{itemize}
\end{corollary}

\begin{corollary}\colabel{5.8}
Let $M$ be a $(T,R)$-bimodule. The functor $G=\Hom({}_TM,{}_T\bullet)$
is separable if and only if $T$ is $M$-separable over $R$, in the sense
of \deref{2.1}. If $M$ is finitely generated and projective as
a left $R$-module, with finite dual basis $\{h_k, p_k\}$, then this
is also equivalent to the existence of $\phi:\ M^*\to {}^*M$
in $W_2$ such that
$$\sum_k (p_k)\phi(h_k)=1$$
\end{corollary}

The Frobenius analog of \coref{5.8} is that ${}_TM_R$ is a Frobenius
bimodule if and only if $M\ot_R\bullet$ is a Frobenius functor. This
was stated explicitely in \cite{CGN99}, where it is also proved that
any additive Frobenius functor between module categories is of this
type. Let us show how this result can be deduced easily from
Propositions \prref{5.5} and \prref{5.6}.

\begin{proposition}\prlabel{5.9}
Let $M$ be a $(T,R)$-bimodule, and consider the functors $F=M\ot_R\bullet$
and $G=\Hom({}_TM,{}_T\bullet)$. The following assertions are equivalent.
\begin{enumerate}
\item $(F,G)$ is a Frobenius pair;
\item ${}_TM$ is finitely generated and projective, and there exist
$e=\sum_i m_i\ot_R f_i\in W_1$ and $\ol{\nu}\in V_1$ such that
\begin{eqnarray}
m&=& \sum_i m_i\ol{\nu}(f_i\cdot m)\eqlabel{5.9.1}\\
f&=& \sum_i\ol{\nu}(f\cdot m_i)f_i\eqlabel{5.9.2}
\end{eqnarray}
for all $m\in M$ and $f\in {}^*M$;
\item ${}_TM_R$ is Frobenius in the sense of \deref{3.8}.
\end{enumerate}
\end{proposition}

\begin{proof}
$1)\Rightarrow 2)$. If $(F,G)$ is Frobenius, then $G=\Hom({}_TM,\bullet)$
is a left
adjoint, and therefore right exact and preserving direct limits,
so ${}_TM$ is necessarily finitely generated and projective.
Let $\zeta$ and $\nu$ be the counit and unit of the adjunction $(G,F)$, i.e.
\begin{equation}\eqlabel{5.9.3}
F(\nu_P)\circ\zeta_{F(P)}=I_{F(P)}~~{\rm and}~~
\nu_{G(Q)}\circ G(\zeta_Q)=I_{G(Q)}
\end{equation}
for all $P\in {}_R{\mathcal M}$ and $Q\in {}_T{\mathcal M}$.
Let $\ol{\nu}=\nu_R\in V_1$ and $e=\sum_i m_i\ot_R f_i=\zeta_T(1)
\in W_1$. Putting $P=R$ in \equref{5.9.3}, we find \equref{5.9.1}.
Then take $Q=T$ in \equref{5.9.3}. Making the identification
$G'=G$ (${}_TM$ is finitely generated projective), we find
$$G(\zeta_T):\ {}^*M\to {}^*M\ot_T M\ot_R {}^*M,~~
G(\zeta_T)(f)=\sum_i f\ot_T m_i\ot_R f_i$$
and $\nu_{G(T)}(G(\zeta_T)(f))=\sum_i\ol{\nu}(f\cdot m_i)f_i$,
so \equref{5.9.2} follows.\\
$2)\Rightarrow 1)$. Let $\nu= \alpha^{-1}(\ol{\nu})$ and
$\zeta=\beta^{-1}(e)$. $\nu$ and $\zeta$ satisfy \equref{5.9.3},
so $(F,G)$ is Frobenius.\\
$2)\Rightarrow 3)$. \equref{5.9.1} implies that $M_R$ is finitely
generated projective. Let $\ol{\phi}=\alpha_1(\ol{\nu})$ and
$\phi=\beta_1(e)$. We easily compute that $\ol{\phi}$ and
$\phi$ are each others inverses.\\
$3)\Rightarrow 2)$. If ${}^*M$ and $M^*$ are isomorphic as $(R,T)$-bimodules,
then there exist $\ol{\phi}\in V_2$ and $\phi\in W_2$ that are each
others inverses. Put $\ol{\nu}=\alpha_1^{-1}(\ol{\phi})$ and
$e=\beta_1^{-1}(\phi)$. Straightforward computations show that
$\ol{\nu}$ and $e$ satisfy (\ref{eq:5.9.1}-\ref{eq:5.9.2}).
\end{proof}

Obviously our results also hold for functors between categories of right
modules. As before, let $M$ be a $(T,R)$-bimodule, and consider the
functors
$$\widetilde{F}=\bullet\ot_T M:\ {\mathcal M}_T\to {\mathcal M}_R~~{\rm and}~~
\widetilde{G}=\Hom(M_R,\bullet_R):\ {\mathcal M}_R\to {\mathcal M}_T$$
Then $(\widetilde{F},\widetilde{G})$ is an adjoint pair. We have a natural
transformation
$$\widetilde{\gamma}:\ \widetilde{G'}=\bullet\ot_R M^*\to G$$
For all $Q\in {\mathcal M}_R$, $\widetilde{\gamma}_Q:\ Q\ot_R M^*\to
\Hom(M_R,Q_R)$ is given by
$$\gamma_Q(q\ot h)(m)=qh(m)$$
We denote $\gamma_Q(q\ot h)=q\cdot h$. The analogs of Propositions
\ref{pr:5.5} and \ref{pr:5.6} are the following:

\begin{proposition}\prlabel{5.10}
With notation as above, we have maps
$$\widetilde{V}=\Nat(\widetilde{G}\widetilde{F},1_{{\mathcal
M}_T})\rTo{\widetilde{\alpha}}
\widetilde{V}_1=\Hom({}_T\Hom(M_R,M_R)_T,{}_TT_T)\rTo{\widetilde{\alpha}_1}W_2$$
and
$$\widetilde{W}=\Nat(1_{{\mathcal
M}_R},\widetilde{F}\widetilde{G})\rTo{\widetilde{\beta}}
\widetilde{W}_1\rTo{\widetilde{\beta}_1}V_2$$
where $\widetilde{W}_1=\{e\in M^*\ot_T M~|~re=er~{\rm for~all~} r\in R\}$.
$\widetilde{\beta}$ is always an isomorphism, $\widetilde{\alpha}$ and
$\widetilde{\alpha}_1$ are isomorphisms if $M_R$ is finitely generated
projective,
and $\widetilde{\beta}_1$ is an isomorphism if ${}_TM$ is finitely generated
projective.
\end{proposition}

\begin{proof}
Completely similar to the proof of Propositions \ref{pr:5.5} and \ref{pr:5.6}.
Let us mention that
$$\alpha_1(\ol{\nu})=\phi~~{\rm with}~~(m)\phi(h)=\ol{\nu}(\gamma_M(m\ot h))$$
$$\beta_1(\sum_i k_i\ot m_i)=\ol{\phi}~~{\rm with}~~\ol{\phi}(f)=\sum_i
k_i((m_i)f) \qed $$
\renewcommand{\qed}{}\end{proof}

As a consequence, we obtain relations between the separability and Frobenius
properties of $F$, $G$, $\widetilde{F}$ and $\widetilde{G}$.

\begin{corollary}\colabel{5.11}
Let $M$ be a $(T,R)$-bimodule, and assume that $M_R$ is finitely generated
projective. Then the following assertions are equivalent:
\begin{itemize}
\item $\widetilde{F}$ is separable;\\
\item there exists $\ol{\nu}\in \widetilde{V}_1$: $\ol{\nu}(I_M)=1_T$;\\
\item there exists $\phi\in W_2$: $\sum_k (p_k)\phi(h_k)=1$;\\
\item $G$ is separable.
\end{itemize}
Here $\{p_k,h_k\}$ is a finite dual basis of $M$ as a right $R$-module.
\end{corollary}

\begin{proof}
The equivalence of the first three statements is obtained in exactly the
same way as \coref{5.7}. The equivalence of the third and the fourth statement
is one of the equivalences in \coref{5.8}.
\end{proof}

\begin{corollary}\colabel{5.12}
Let $M$ be a $(T,R)$-bimodule. The following statements are equivalent:
\begin{itemize}
\item $\widetilde{G}$ is separable;
\item there exists $e=\sum_i k_i\ot_T m_i\in \widetilde{W}_1$ such that
$\sum_i k_i(m_i)=1_R$.
\end{itemize}
If ${}_TM$ is finitely generated projective, then they are also equivalent to
\begin{itemize}
\item There exists $\ol{\phi}\in V_2$ such that $\sum_j\ol{\phi}(g_j)n_j=1_R$;
\item $F$ is separable.
\end{itemize}
Here $\{n_j,g_j\}$ is a finite dual basis of $M$ as a left $T$-module.
\end{corollary}

\begin{corollary}\colabel{5.13}
Let $M$ be a $(T,R)$-bimodule. The following statements are equivalent:
\begin{itemize}
\item $(\widetilde{F},\widetilde{G})$ is a Frobenius pair and
$M_R$ is finitely generated projective;
\item $M_R$ is finitely generated projective and there exist $e=\sum_i
k_i\ot_T m_i\in
\widetilde{W}_1$ and $\ol{\nu}\in \widetilde{V}_1$ such that
\begin{equation}\eqlabel{5.13.1}
m=\sum_i \ol{\nu}(q\cdot k_i)m_i~~{\rm and}~~
h=\sum_i k_i\ol{\nu}(m_i\cdot f)
\end{equation}
for all $m\in M$ and $f\in M^*$;
\item ${}_TM_R$ is a Frobenius bimodule.
\end{itemize}
\end{corollary}

We now address the following problem: assume that we know that
${}_TM_R$ is a Frobenius bimodule. When is $T$ $M$-separable
over $R$?  In view of the above considerations, it suffices
to apply \prref{5.3}. We first need a Lemma.

\begin{lemma}\lelabel{5.14}
Let $M$ be a $(T,R)$-bimodule, and $F=M\ot_R\bullet$. Then
$$\Nat(F,F)\cong \Hom({}_TM_R,{}_TM_R)$$
\end{lemma}

\begin{proof}
Let $\alpha:\ F\to F$ be a natural transformation. Then $\alpha_R:\ M\to M$
is left $T$-linear, and right $R$-linear since $\alpha$ is natural.
Given a $(T,R)$-linear map $\alpha_R:\ M\to M$, we define a natural
transformation $\alpha:\ F\to F$ as follows: $\alpha_P=\alpha_R\ot I_P$.
\end{proof}

\begin{proposition}\prlabel{5.15}
Assume that $M$ is a Frobenius $(T,R)$-bimodule, and let
$e=\sum_i m_i\ot_R f_i\in W_1$ and $\ol{\nu}\in V_1$ be as in the second
statement of \prref{5.9}. Then $T$ is $M$-separable over $R$ (i.e. $G$
is separable) if and only if there exists a $(T,R)$-linear map
$\ol{\alpha}:\ M\to M$ such that
$$\sum_i (\ol{\alpha}(m_i))f_i=1_T$$
\end{proposition}

This proposition recovers \cite[Theorem 2.15]{K}.  One can similarly
recover \cite[Theorem 4.1]{K99}.

\end{section}
%%%%%%%%%%%%%%%%%%%%%%%%%%%%%%%%%%%%%%%%%%%%%%%%%%%%%%%%%%
%%%%%%%%%%%%   BIBLIOGRAPHY  %%%%%%%%%%%%%%%%%%%%%%%%%%%%%
%%%%%%%%%%%%%%%%%%%%%%%%%%%%%%%%%%%%%%%%%%%%%%%%%%%%%%%%%%

\bibliographystyle{amsalpha}

\end{document}